\documentclass[11pt]{article}

\newcommand{\bc}{\begin{center}}
\newcommand{\ec}{\end{center}}
\newcommand{\be}{\begin{equation}}
\newcommand{\ee}{\end{equation}}
\newcommand{\ba}{\begin{array}}
\newcommand{\ea}{\end{array}}
\newcommand{\bea}{\begin{eqnarray}}
\newcommand{\eea}{\end{eqnarray}}
\newcommand{\edc}{\end{document}}

\def\f{\varphi}
\def\p{\psi}
\def\s{\sigma}

\def\O{\Omega}
\def\w{\omega}

\def\s{\sigma}

\def\l{\lambda}

\def\a{\alpha}
\def\w{\omega}
\def\O{\Omega}

\def\b{\beta}
\def\P{\Phi}

\def \d {\partial}
\def \L {\Lambda}

\def \a {\alpha}

\begin{document}
\thispagestyle{empty}
\begin{center}

{\bf ON NON-VOLTERRA QUADRATIC STOCHASTIC OPERATORS\\ GENERATED BY
A PRODUCT MEASURE}\\
\vspace{1.3cm}

U.A. Rozikov\footnote{rozikovu@yandex.ru}\\
{\it Institute of Mathematics, 29, F. Hodjaev str., Tashkent, 700125, Uzbekistan\\
and \\
The Abdus Salam International Centre for Theoretical Physics, Triesty, Italy.}\\
and\\[1em]
N.B. Shamsiddinov \footnote{shnb$_-$siz@yahoo.com}\\{\it Institute
of Irrigation and Agricultural Mechanization
Engineers,\\ Tashkent, 700052, Uzbekistan.}\\
\end{center}
 \vspace{1cm} \centerline{\bf Abstract}
\bigskip
\baselineskip=18pt

In this paper we describe a wide class of non-Volterra quadratic
stochastic operators using  N. Ganikhadjaev's construction of
quadratic stochastic operators. By the construction these
operators depend on a probability measure $\mu$ being defined on
the  set of all configurations which are given on a graph $G.$  We
show that if $\mu$ is  the product of probability measures being
defined on each maximal connected subgraphs of $G$ then
corresponding non-Volterra operator can be reduced to $m$ number
(where $m$ is the number of maximal connected subgraphs of $G$) of
Volterra operators defined on the maximal connected subgraphs. Our
result allows to study a wide class of non-Volterra operators in
the framework of the well known theory of Volterra quadratic
stochastic operators.


\section{Introduction}

It is well known that the principles of biological inheritance,
initiated by Mendel allow  an exact mathematical formulation. For
this reason classical genetics can be regarded as a mathematical
discipline. The fundamental investigations on these problems were
carried out by S.N.Bernstein [2] and Reiersol [13].

The main mathematical apparatus of such investigations, to our
knowledge, is the theory of quadratic stochastic operators.
 Such operator frequently arise in many models of mathematical
genetics [1-11], [13], [15],[16].

The quadratic stochastic operator (QSO) is a mapping of the
simplex
$$ S^{n-1}=\{x=(x_1,...,x_n)\in R^n: x_i\geq 0, \sum^n_{i=1}x_i=1 \} $$
into itself, of the form
$$ V: x_k'=\sum^n_{i,j=1}p_{ij,k}x_ix_j, \ \ (k=1,...,n), \eqno(1)$$
where $p_{ij,k}$ are coefficients of heredity and
$$ p_{ij,k}\geq 0, \ \ \sum^n_{k=1}p_{ij,k}=1, \ \ (i,j,k=1,...,n). \eqno(2)$$
Note that each element $x\in S^{n-1}$ is a probability
distribution on $E=\{1,...,n\}.$

The population evolves by starting from an arbitrary state
(probability distribution on $E$) $x\in S^{n-1}$ then passing to
the state $Vx$ (in the next ``generation''), then to the state
$V^2x$, and so on.

For a given $x^{(0)}\in S^{n-1}$ the trajectory $\{x^{(l)}\},
l=0,1,2,...$ of QSO (1) is defined by $x^{(l+1)}=V(x^{(l)}),$
where $l=0,1,2,...$  One of the main problems in mathematical
biology consists in the study of the asymptotical behavior of the
trajectories. This problem was fully solved [1-7] for Volterra QSO
which is defined by (1), (2) and the additional assumption
$$ p_{ij,k}=0, \ \ {\rm if} \ \ k\notin \{i,j\}. \eqno(3)$$
The biological treatment of condition (3) is rather clear: the
offspring repeats the genotype of one of its parents.

In paper [6] the general form of Volterra QSO $V:
x=(x_1,...,x_n)\in S^{n-1}\to V(x)=x'=(x'_1,...,x'_n)\in S^{n-1}$
is given:

$$x'_k=x_k\bigg(1+\sum_{i=1}^na_{ki}x_i\bigg), \eqno(4)$$
where $a_{ki}=2p_{ik,k}-1$ for $i\ne k$ and $a_{kk}=0.$ Moreover
$a_{ki}=-a_{ik}$ and $|a_{ki}|\leq 1.$

In papers  [6], [7] the theory of QSO (4) was  developed using
theory of the Lyapunov functions and tournaments. But non-Volterra
QSOs (i.e. which do not satisfy the condition (3)) were not in
completely studied. Because there is no any general theory which
can be applied for investigation of non-Volterra operators. To the
best of our knowledge, there are few papers devoted to such
operators (see e.g. [1], [4]).

In papers [3],[5] a constructive description of QSOs is given.
This construction depends on a probability measure $\mu$ and
cardinality of a set of cells (configurations) which  can be
finite or continual. In this paper we describe QSOs using the
construction of QSO for the general finite graph and probability
measure $\mu$ which is product of measures defined on maximal
subgraphs of the graph. We show that if $\mu$ is given by the
product of the probability measures then corresponding
non-Volterra operator can be reduced to $m$ number (where $m$ is
the number of maximal connected subgraphs) of Volterra operators
defined on the maximal connected subgraphs. Then using the well
known theory of Volterra operators we study the asymptotical
behavior of the trajectories of the non-Volterra operators.

\section {Construction of QSO}

Each quadratic operator $V$ can be uniquely defined by a cubic
matrix ${\mathbf P}\equiv {\mathbf P}(V)=\{p_{ij,k}\}^n_{i,j,k=1}$
with condition (2). Usually [1], [2], [4], [6-11] the matrix
${\mathbf P}$ is known. In [3], [5]  a constructive description of
${\mathbf P}$ is given.

Now we recall this construction.

Let $G=(\L, L)$ be a finite graph without loops and multiple
edges, where $\L$ is the set of vertexes and $L$ is the set of
edges of the graph.

Furthermore, let $\Phi$ be a finite set, called the set of alleles
(in problems of statistical mechanics, $\P$ is called the  range
of spin). The function $\s:\L\to\P$ is called a cell (in mechanics
it is called configuration). Denote by $\O$ the set of all cells,
this set corresponds to $E$ in (1). Let $S(\L, \P)$ be the set of
all probability measures defined on the finite set $\O.$

Let $\{\L_i, i=1,...,m\}$ be the set of maximal connected
subgraphs (components) of the graph $G.$ For  a configuration
$\s\in \O$ denote by  $\s(M)$ its "projection" (or "restruction")
to $M\subset \L$ : $\s(M) =\{\s(x)\}_{x\in M}.$ Fix two cells
$\s_1, \s_2\in \O,$ and put
$$\O(G, \s_1,\s_2)=\{\s\in \O: \s(\L_i)=\s_1(\L_i) \ \ {\rm or} \ \
\s(\L_i)=\s_2(\L_i) \ \ {\rm for \ \ all} \ \ i=1,...,m\}.$$ Now
let $\mu\in S(\L,\P)$ be a probability measure defined on $\O$
such that $\mu(\s)>0$ for any cell $\s\in \O;$ i.e $\mu$ is a
Gibbs measure with some potential [12], [14]. The heredity
coefficients $p_{\s_1\s_2,\s}$ are defined as
$$p_{\s_1\s_2,\s}=\left\{\begin{array}{ll}
{\mu(\s)\over \mu(\O(G,\s_1,\s_2))}, \ \ {\rm if} \ \ \s\in
\O(G,\s_1,\s_2),\\[3mm]
0 \ \ {\rm otherwise}.\\
\end{array}\right.\eqno (5)$$
Obviously, $p_{\s_1\s_2,\s}\geq 0,$
$p_{\s_1\s_2,\s}=p_{\s_2\s_1,\s}$ and $\sum_{\s\in \O}
p_{\s_1\s_2,\s}=1$ for all $\s_1,\s_2\in \O.$

The QSO $V\equiv V_\mu$ acting on the simplex $S(\L,\P)$ and
determined by coefficients (5) is defined as follows: for an
arbitrary measure $\l\in S(\L,\P)$, the measure $V(\l)=\l'\in
S(\L,\P)$ is defined by the equality
$$ \l'(\s)=
\sum_{\s_1,\s_2\in \O} p_{\s_1\s_2,\s} \l(\s_1)\l(\s_2) \eqno(6)$$
for any cell $\s\in \O.$

The QSO construction is also closely related to the graph
structure on the set $\L.$

\vskip 0.4 truecm {\bf Theorem 1.} [3] {\it The QSO (5) is
Volterra if and only if the graph $G$ is connected.} \vskip 0.3
truecm

\section{Non-Volterra  QSO}

In this section we describe a condition on measure $\mu$ under
which the QSO $V_{\mu}$ generated by $\mu$ can be studied using
the theory of Volterra QSO.

Let $G=(\L,L)$ be a finite graph and $\{\L_i, i=1,...,m\}$ the set
of all maximal connected subgraphs of $G$. Denote by
$\O_i=\P^{\L_i}$ the set of all configurations defined on $\L_i,$
$i=1,...,m.$ Let $\mu_i$ be a probability measure defined on
$\O_i,$ such that $\mu_i(\s)>0$ for any $\s\in \O_i,$ $i=1,...,m.$

Consider probability measure $\mu$ on $\O=\O_1\times\dots\times
\O_m$ defined as

$$ \mu(\s)=\prod^m_{i=1}\mu_i(\s_i), \eqno (7)$$
where $\s=(\s_1,...,\s_m),$ with $\s_i\in \O_i, i=1,...,m.$

By Theorem 1 if $m=1$ then QSO constructed on $G$ is Volterra QSO.
\vskip 0.4 truecm

{\bf Theorem 2.} {\it The QSO constructed by (5) with measure (7)
is reducible to $m$ separate Volterra QSOs.}

\vskip 0.4 truecm

{\bf Proof.} Take $\f=(\f_1,...,\f_m), \p=(\p_1,...,\p_m)\in \O.$
By construction
$$\O(G, \f,\p)=\bigg\{\s=(\s_1,..,\s_m)\in \O: \s_i\in \{\f_i,\p_i\},
i=1,...,m\bigg\}$$ and
$$p_{\f\p,\s}=\left\{\begin{array}{ll}
\prod^m_{i=1}{\mu_i(\s_i)\over \mu_i(\f_i)+\mu_i(\p_i)}, \ \ {\rm
if} \ \ \s\in
\O(G,\f,\p),\\[3mm]
0 \ \ {\rm otherwise},\\
\end{array}\right.\eqno (8)$$
where we used the following equality
$$\mu(\O(G,\f,\p))=\sum_{\s_1,...,\s_m:\atop \s_i\in \{\f_i,
\p_i\},
i=1,...,m}\prod^m_{i=1}\mu_i(\s_i)=\prod_{i=1}^m\big(\mu_i(\f_i)+\mu_i(\p_i)\big).$$
Thus QSO (6) generated by measure (7) can be written as
$$\l'(\s)=\l'(\s_1,...,\s_m)=\sum_{\f=(\f_1,...,\f_m): \f_i\in \O_i\atop \p=(\p_1,...,\p_m): \p_i\in \O_i}
\prod^m_{i=1}{\mu_i(\s_i){\mathbf 1}_{(\s_i\in
\{\f_i,\p_i\})}\over
\mu_i(\f_i)+\mu_i(\p_i)}\l(\f)\l(\p).\eqno(9)$$ Denote
$$X_{i,\w}=\sum_{\s\in \O:\atop
\s_i=\w}\l(\s)=\sum_{\s_1,...,\s_{i-1},\s_{i+1},...,\s_m\atop
\s_k\in \O_k, k\ne i}\l(\s_1,...,\s_{i-1},\w,\s_{i+1},...,\s_m).
\eqno(10)$$ From (9) we have
$$X'_{i,\w}=\sum_{\s\in \O:\atop \s_i=\w}\l'(\s)=\sum_{\s\in \O:\atop
\s_i=\w}\bigg[\sum_{\f_1,...,\f_{i-1},\f_{i+1},...,\f_m\atop \p_1,
..., \p_m}{\mu_i(\w)\over \mu_i(\w)+\mu_i(\p_i)}\times
$$
$$
\prod_{j=1\atop j\ne i}^m{\mu_j(\s_j){\mathbf 1}_{(\s_j\in
\{\f_j,\p_j\})}\over
\mu_j(\f_j)+\mu_j(\p_j)}\l(\f_1,...,\f_{i-1},\w,\f_{i+1},...,\f_m)\l(\p_1,...,\p_m)+$$
$$
\sum_{\f_1,...,\f_m\atop \p_1, ...,\p_{i-1},\p_{i+1},...,
\p_m}{\mu_i(\w)\over \mu_i(\f_i)+\mu_i(\w)}\times$$
$$
\prod_{j=1\atop j\ne i}^m{\mu_j(\s_j){\mathbf 1}_{(\s_j\in
\{\f_j,\p_j\})}\over
\mu_j(\f_j)+\mu_j(\p_j)}\l(\f_1,...,\f_m)\l(\p_1,...,\p_{i-1},\w,\p_{i+1},...,\p_m)\bigg]=$$
$$ 2\sum_{\f_1,...,\f_{i-1},\f_{i+1},...,\f_m\atop
\p_1, ..., \p_m}{\mu_i(\w)\over \mu_i(\w)+\mu_i(\p_i)}\times
$$
$$\sum_{\s_1,...,\s_{i-1},\s_{i+1},...,\s_m}
\prod_{j=1\atop j\ne i}^m{\mu_j(\s_j){\mathbf 1}_{(\s_j\in
\{\f_j,\p_j\})}\over
\mu_j(\f_j)+\mu_j(\p_j)}\l(\f_1,...,\f_{i-1},\w,\f_{i+1},...,\f_m)\l(\p_1,...,\p_m).\eqno(11)$$

Note that
$$\sum_{\s_1,...,\s_{i-1},\s_{i+1},...,\s_m}
\prod_{j=1\atop j\ne i}^m{\mu_j(\s_j){\mathbf 1}_{(\s_j\in
\{\f_j,\p_j\})}\over \mu_j(\f_j)+\mu_j(\p_j)}=1.$$

Thus from (11) we have
$${\rm RHS\ \ of\ \ (11)}=$$
$$ 2\sum_{\f_1,...,\f_{i-1},\f_{i+1},...,\f_m\atop \p_1, ...,
\p_m}{\mu_i(\w)\over \mu_i(\w)+\mu_i(\p_i)}
\l(\f_1,...,\f_{i-1},\w,\f_{i+1},...,\f_m)\l(\p_1,...,\p_m)=$$
$$\sum_{\f_1,...,\f_{i-1},\f_{i+1},...,\f_m\atop \p_1,
...,\p_{i-1},\p_{i+1},...,\p_m}\l(\f_1,...,\f_{i-1},\w,\f_{i+1},...,\f_m)
\l(\p_1,...,\p_{i-1},\w,\p_{i+1},...,\p_m)+$$
$$2\sum_{\p_i\in \O_i\setminus \w}{\mu_i(\w)\over
\mu_i(\w)+\mu_i(\p_i)}
\sum_{\f_1,...,\f_{i-1},\f_{i+1},...,\f_m\atop \p_1,
...,\p_{i-1},\p_{i+1},..., \p_m}
\l(\f_1,...,\f_{i-1},\w,\f_{i+1},...,\f_m)\l(\p_1,...,\p_m)=$$
$$X_{i,\w}^2+\sum_{\p\in \O_i\setminus \w}{2\mu_i(\w)\over
\mu_i(\w)+\mu_i(\p)}X_{i,\w}X_{i,\p}.$$

Thus operator (9) can be rewritten as
$$X'_{i,\w}=X_{i,\w}\bigg(X_{i,\w}+\sum_{\p\in \O_i\setminus \w}{2\mu_i(\w)\over
\mu_i(\w)+\mu_i(\p)}X_{i,\p}\bigg), \eqno(12)$$ where $X_{i,\w}$
is defined by (10), $\w\in \O_i, i=1,...,m.$

Note that $\sum_{\w\in \O_i}X_{i,\w}=1$ for any $i=1,...,m.$ Using
this equality from (12) we obtain
$$X'_{i,\w}=X_{i,\w}\bigg(1+\sum_{\p\in \O_i}{\mu_i(\w)-\mu_i(\p)\over
\mu_i(\w)+\mu_i(\p)}X_{i,\p}\bigg). \eqno(13)$$ Comparing this
with (4) one can see that for each fixed $i$ ($i=1,...,m$) the
operator (13) is Volterra operator $V^{(i)}:S^{|\O_i|-1}\to
S^{|\O_i|-1}.$ The theorem is proved.

\vskip 0.4 truecm

{\bf Remark.} The set $\P$ which we used for the set of spin
values on $G$ is the same for any configuration on $\O_i$,
$i=1,...,m.$ Note that Theorem 2 can be easily extended for more
general case: when each subgraph $\L_i$ has its own $\Phi_i.$

\section{The behavior of the trajectories}

In this section using Volterra QSOs (13) we shall describe the
behavior of trajectories of non-Volterra QSO (9).

Denote $a^{(i)}_{\f,\p}={\mu_i(\w)-\mu_i(\p)\over
\mu_i(\w)+\mu_i(\p)}.$ It is easy to see that for each fixed $i\in
\{1,...,m\}$ the coefficients  $a^{(i)}_{\f,\p}$ satisfy the
following properties
$$a^{(i)}_{\f,\p}=-a^{(i)}_{\p,\f}, \ \ |a^{(i)}_{\f,\p}|\leq 1.
\eqno (14)$$ Thus QSO (13) has the same form with (4).

If for any $i\in \{1,...,m\}$ the asymptotical behavior of
trajectories of QSO $V^{(i)}$ (i.e. (13)) is known, say
$X^{(l)}_{i,\w}\to X^*_{i,\w},$ $l\to \infty$, then asymptotical
behavior of $V$ (i.e. (9)), say $\l^{(l)}(\s)\to \l^*(\s),$ $l\to
\infty$, can be found from system of linear equations
$$\sum_{\s\in \O: \s_i=\w}\l^*(\s)=X^*_{i,\w}, \ \ \w\in \O_i,
i=1,...,m. \eqno(15)$$

Now we recall known results ([6],[7]) about trajectories of (4)
(i.e. (13) for fixed $i$).

Put
$$ \d S^{n-1}=\{x=(x_1,...,x_n)\in S^{n-1}: \prod^n_{j=1}x_j=0\},
\ \ {\rm int}S^{n-1}=S^{n-1}\setminus \d S^{n-1}.$$ Let
$x^{(0)}\in S^{n-1}$ be the initial point. Denote by $\nu(x^0)$
the set of limit points of the trajectory $\{x^{(l)}\}$ of the
operator (4).

\vskip 0.4 truecm

{\bf Theorem 3.} [6]. 1) {\it If $x^{(0)}\in {\rm int}S^{n-1}$ is
not a fixed point (i.e. $Vx^{(0)}\ne x^{(0)}$), then}
$\nu(x^{(0)})\subset \d S^{n-1}.$

2) {\it The set $\nu(x^{(0)})$ either consists of a single point
or is infinite.}

3) {\it If QSO (4) has an isolated fixed point $x^*\in {\rm
int}S^{n-1},$ then for any initial point $x^{(0)}\ne x^*$, the
trajectory $\{x^{(l)}\}$ does not converge.}

\vskip 0.4 truecm

{\bf Corollary 4.} {\it The set of limit points  ${\tilde
\nu}(\l^0)$ of the QSO (9) has properties 1)-3) mentioned in
Theorem 3.}

\vskip 0.4 truecm

{\bf Proof.} 1) By Theorem 3 there is at least one ${\tilde \w}\in
\O_i$ with $X^*_{i,{\tilde \w}}=0.$ Since $\l(\s)\geq 0$ from (15)
we get that $\l^*(\s)=0$ for all $\s$ such that $\s_i={\tilde
\w}.$ This completes the proof of the property 1). Properties
2),3) also follow from (15).

Now we shall recall some notions from theory of tournaments (see
[6]) for operator (4). Assume $a_{ki}\ne 0$ for $k\ne i.$ Along
with (4) consider the complete graph $G_n$  with $n$ vertices.
Specify a direction on the edges of $G_n$ as follows: the edge
joining vertices $k$ and $i$ is directed from the $k$th to the
$i$th vertex if $a_{ki}<0,$ and has the opposite direction if
$a_{ki}>0.$ The directed graph thus obtained is called a
tournament and denoted by $T_n.$  A tournament is said to be
strong if it is possible to go from any vertex to any other vertex
with directions taken into account.

Further, after a suitable renumbering of the vertices of $T_n$ we
can assume that the subtournament $T_r$ contains the first $r$
vertices of $T_n$ as its vertices. Obviously, $r\leq n,$ and $r=n$
if and only if $T_n$ is a strong tournament.

\vskip 0.4 truecm

{\bf Theorem 5.}[6]. {\it Suppose that $T_n$ is not strong, and
let $x^{(0)}\in {\rm int}S^{n-1}.$ If $j>r,$ then $x^{(l)}_j\to
0,$ at the rate of a geometric progression as $l\to \infty.$}

As corollary of Theorem 5 we have

\vskip 0.4 truecm

 {\bf Corollary 6.} {\it Suppose that there is a
$i_0\in \{1,...,m\}$ such that the corresponding tournament
$T^{(i_0)}$ of the operator $V^{(i_0)}$ is not strong, and assume
$\l^0\in {\rm int}S^{n-1}.$ Then there is a subset ${\tilde
\O}_{i_0}\subset \O_{i_0}$ such that $\l^{(l)}(\s)\to 0,$ for
$\s\in \O$ with $\s_{i_0}\in {\tilde \O}_{i_0},$ at the rate of a
geometric progression as $l\to \infty.$}

\vskip 0.4 truecm {\bf Example.} Consider graph $G=(\L,L)$ with
$\L=\{1,2\}$ and $L=\emptyset.$ Take $\P=\{A, a\}.$ Then
non-Volterra QSO (9)  has the form
$$\begin{array}{llll}
x_1'=x_1^2+2\b_1x_1x_2+2\a_1x_1x_3+2\a_1\b_1x_1x_4+2\a_1\b_1x_2x_3\\[2mm]
x_2'=x_2^2+2\b_2x_1x_2+2\a_1\b_2x_2x_3+2\a_1x_2x_4+2\a_1\b_2x_1x_4\\[2mm]
x_3'=x_3^2+2\a_2x_1x_3+2\a_2\b_1x_2x_3+2\b_1x_3x_4+2\a_2\b_1x_1x_4\\[2mm]
x_4'=x_4^2+2\a_2\b_2x_1x_4+2\a_2x_2x_4+2\b_2x_3x_4+2\a_2\b_2x_2x_3\\[2mm]
\end{array}\eqno(16)$$
where $\mu_1=(\a_1, \a_2), \a_j\geq 0, j=1,2,  \a_1+\a_2=1;
\mu_2=(\b_1,\b_2), \b_j\geq 0, j=1,2,  \b_1+\b_2=1.$

Putting $x_1+x_2=X_{1,1}$, $x_3+x_4=X_{1,2}$ and
$x_1+x_3=X_{2,1},$ $x_2+x_4=X_{2,2}$ we get the Volterra operators
like (13):
$$\begin{array}{ll}
X_{1,1}'=X_{1,1}\bigg(1+(2\a_1-1)X_{1,2}\bigg)\\[2mm]
X_{1,2}'=X_{1,2}\bigg(1+(2\a_2-1)X_{1,1}\bigg)\\[2mm]
\end{array}\eqno(17)$$
and
$$\begin{array}{ll}
X_{2,1}'=X_{2,1}\bigg(1+(2\b_1-1)X_{2,2}\bigg)\\[2mm]
X_{2,2}'=X_{2,2}\bigg(1+(2\b_2-1)X_{2,1}\bigg)\\[2mm]
\end{array}\eqno(18)$$

For this example from Corollary 4 and 6 we obtain

{\bf Corollary 7.} 1. {\it Any trajectory of non-Volterra QSO (16)
has following limit}
$$\lim_{l\to \infty}x^{(l)}= \left\{\begin{array}{llll}
(1,0,0,0), \ \ {\rm if} \ \ 2\a_1>1, 2\b_1>1 ,\\[3mm]
(0,1,0,0), \ \ {\rm if} \ \ 2\a_1>1, 2\b_1<1 ,\\[3mm]
(0,0,1,0), \ \ {\rm if} \ \ 2\a_1<1, 2\b_1>1 ,\\[3mm]
(0,0,0,1), \ \ {\rm if} \ \ 2\a_1<1, 2\b_1<1 .\\[3mm]
\end{array}\right.$$

2. {\it If $2\b_1=1$ then  $S_1=\{x: x_3=x_4=0\}$ and $S_2=\{x:
x_1=x_2=0\}$ are the sets of fixed points for (16) and for any}
$x^{(0)}\notin S_1\cup S_2$
$$\lim_{l\to \infty}x^{(l)}\in \left\{\begin{array}{ll}
S_1, \ \ {\rm if} \ \ 2\a_1>1 ,\\[3mm]
S_2 \ \ {\rm if} \ \ 2\a_1<1.\\
\end{array}\right.$$

3. {\it If $2\a_1=1$ then $S_3=\{x: x_2=x_4=0\}$ and $S_4=\{x:
x_1=x_3=0\}$ are the sets of fixed points for (16) and for any}
$x^{(0)}\notin S_3\cup S_4$
$$\lim_{l\to \infty}x^{(l)}\in \left\{\begin{array}{ll}
S_3, \ \ {\rm if} \ \ 2\b_1>1 ,\\[3mm]
S_4 \ \ {\rm if} \ \ 2\b_1<1.\\
\end{array}\right.$$

4. {\it If $2\a_1=2b_1=1$ then  $S_5=\{x: x_2=x_4, x_1=x_3\}$ and
$S_6=\{x: x_1=x_2, x_3=x_4\}$ are the sets of fixed points for
(16).}
 \vskip 0.4
truecm {\bf Acknowledgments.} The final part of this work was done
at the Abdus Salam International Centre for Theoretical Physics
(ICTP), Trieste, Italy and the first author thank ICTP for
providing financial support and all facilities. The first author
also was particularly supported by NATO Reintegration Grant:
FEL.RIG. 980771.

\vskip 0.4 truecm

{\bf References}

1. Anorov O.U., Shamsiddinov N.B., On the reduction of a class of
quadratic stochastic operators to Volterra operators. {\it Uzbek
Math. Jour.} No.3-4 : 9-12 (2001). (Russian)

2. Bernshtein S.N.,  Solution of a mathematical problem connected
with the theory of heredity, {\it Uch. Zap. Nauchno-Issled. kaf.
Ukr. Otd. Mat.,} {\bf 1} : 83-115 (1924).

3. Ganikhodjaev N.N., An application of the theory of Gibbs
distributions to mathematical genetics, {\it Doklady Math.} {\bf
61}: 321-323 (2000).

4. Ganikhodjaev N.N., Mukhitdinov R.T.,  On a class of
non-Volterra quadratic operators, {\it Uzbek Math. Jour.} No. 3-4:
65-69 (2003).

5. Ganikhodjaev N.N., Rozikov U.A.,  On quadratic stochastic
operators generated by Gibbs distributions, {\it Regular and
Chaotic  Dynamics.} {\bf 11}: No. 3 (2006).

6. Ganikhodzhaev R.N., Quadratic stochastic operators, Lyapunov
functions and tournaments, {\it Russian Acad. Sci. Sbornik Math.}
{\bf 76}: 489-506 (1993).

7. Ganikhodzhaev R.N., A chart of fixed points and Lyapunov
functions for a class of discrete dynamical systems. {\it Math.
Notes} {\bf 56}: 1125-1131 (1994).

8. Elson R.C., Stewart J., A general model for the genetic
analysis of pedigree data, {\it Hum. Hered. } {\bf 21}: 523-542
(1971).

9. Jenks, R.D., Quadratic differential systems for interactive
population models, {\it J. Differ. Equations}, {\bf 5}: 497-514
(1969).

10. Lyubich Yu.I., Basic concepts and theorems of the evolutionary
genetics of free populations, {\it Russian Math. Surveys}, {\bf
26}: 51-123 (1971).

11. Lyubich Yu.I., A topological approach to a problem in
mathematical genetics.{\it Russian Math. Surveys} {\bf 34}: 60-66
(1979).

12. Preston C., {\it Gibbs states on countable sets} (Cambridge
University Press, London 1974).

13. Reiersol O., Genetic algebras studied recursively and by means
of differential operators. {\it Math. Scand.} {\bf 10}: 25-44
(1962).

14. Sinai Ya. G., {\it Theory of phase transitions: Rigorous
Results} (Pergamon, Oxford, 1982).

15. Ulam S.M., {\it Problems in Modern Math.}, New York; Wiley,
1964.

16. Zakharevich M.I., The behavior of trajectories and the ergodic
hypothesis for quadratic mappings of a simplex. {\it Russian Math.
Surveys}, {\bf 33}: 207-208 (1978).

\end{document}